\documentclass[a4paper, 11pt]{amsart}

\usepackage{amsfonts} 
\usepackage{amsmath} 
\usepackage{amsthm} 
\usepackage{amssymb} 
\usepackage{enumerate} 
\usepackage{tikz-cd} 

\title[On the positivity of Chern--Weil forms]{A note on Griffiths' conjecture about the positivity of Chern--Weil forms}
\author{Filippo Fagioli}
\address{Filippo Fagioli \\ Dipartimento di Matematica \lq\lq Guido Castelnuovo\rq\rq{} \\ Sapienza Universit\`a di Roma \\ Piazzale Aldo Moro 5 \\ I-00185 Roma.}
\email{fagioli@mat.uniroma1.it} 
\date{January 1, 2022}

\theoremstyle{plain}
\newtheorem{thm}{Theorem}[section]
\newtheorem{cor}[thm]{Corollary}
\newtheorem{lem}[thm]{Lemma}
\newtheorem{prop}[thm]{Proposition}
\newtheorem{quest}{Question}

\theoremstyle{remark}
\newtheorem{rem}[thm]{Remark}
\newtheorem{ex}[thm]{Example}

\theoremstyle{definition}
\newtheorem{defn}[thm]{Definition}

\newcommand{\N}{\mathbb{N}}
\newcommand{\Z}{\mathbb{Z}}
\newcommand{\Q}{\mathbb{Q}}
\newcommand{\R}{\mathbb{R}}
\newcommand{\C}{\mathbb{C}}
\newcommand{\PP}{\mathbb{P}}
\newcommand{\GG}{\mathbb{G}}
\newcommand{\FF}{\mathbb{F}}
\newcommand{\OO}{\mathcal{O}}
\newcommand{\pis}{\pi_{*}}

\begin{document}
\bibliographystyle{amsalpha}
	
	\begin{abstract}
	Let $ (E,h) $ be a Griffiths semipositive Hermitian holomorphic vector bundle of rank $ 3 $ over a complex manifold.
	In this paper, we prove the positivity of the characteristic differential form $ c_1(E,h) \wedge c_2(E,h) - c_3(E,h) $, thus providing a new evidence towards a conjecture by Griffiths about the positivity of the Schur polynomials in the Chern forms of Griffiths semipositive vector bundles.
	As a consequence, we establish a new chain of inequalities between Chern forms.
	Moreover, we point out how to obtain the positivity of the second Chern form $ c_2(E,h) $ in any rank, starting from the well-known positivity of such form if $ (E,h) $ is just Griffiths positive of rank $ 2 $.
	The final part of the paper gives an overview on the state of the art of Griffiths' conjecture, collecting several remarks and open questions.
	\end{abstract}

\maketitle

	\section*{Introduction}
	Let $ X $ be a complex manifold of dimension $ n $ and let $ (E,h) $ be a Hermitian holomorphic vector bundle of rank $ r $ over $ X $.
	Starting from the Chern curvature tensor $ \Theta(E,h) \in \mathcal{A}^{1,1}\bigl(X,\operatorname{End}(E)\bigr) $, let us consider the corresponding Chern forms on $ X $ defined, for $ 0 \le k \le r $, as
	\[
	c_k(E,h) = \operatorname{tr}_{\operatorname{End}({\Lambda}^kE)}\left( {\Lambda}^k \frac{i}{2\pi} \Theta(E,h) \right).
	\]
	By the Chern--Weil theory the form $ c_k(E,h) $ represents the Chern class $ c_k(E) $ of the vector bundle $ E $.
	
	Recall that $ (E,h) \to X $ is called \emph{Griffiths semipositive} (resp. \emph{Griffiths positive}) if for every $x\in X$, $v\in E_x$, $\tau\in T_{X,x}$ we have
	$$
	\langle\Theta(E,h)_x\cdot v,v\rangle_h(\tau,\bar\tau)\ge 0
	$$
	(resp. $ > 0 $ and $=0$ if and only if $v$ or $\tau$ is the zero vector), and that a holomorphic vector bundle $ \mathcal{E} \to X $ is called \emph{ample} if the tautological line bundle $ \OO_{\PP(\mathcal{E}^{\vee})}(1) \to \PP(\mathcal{E}^{\vee}) $ is ample. The notation $ \PP(\mathcal{E}^{\vee}) $ stands for the projective bundle of hyperplanes in $ \mathcal{E} $.
	It is well known that a Griffiths positive vector bundle over a compact complex manifold is ample (the converse is not known in general, but it is a conjecture) and that a globally generated vector bundle (\textsl{i.e}, a bundle whose fibers are spanned by the global holomorphic sections) can be equipped with a Hermitian metric which makes it Griffiths semipositive.
	
	It is natural to expect that certain conditions of positivity on the vector bundle (such as ampleness or Griffiths positivity) impose, in turn, the positivity of objects that derive from it.
	This kind of positivity issues can be placed at both algebraic and differential geometric level.
	In both ways, they were considered by Griffiths in \cite{Gri69}, where it was raised the question of determining the characteristic classes/forms built out of an ample/Griffiths positive vector bundle that are positive in the appropriate sense.
	We refer to \cite{Li20,Xia20,Fin20,DF20} as the most recent works offering a detailed explanation of the results regarding Griffiths' program over the years.
	
	At the algebraic level, one of the main results is the celebrated Fulton--Lazarsfeld theorem \cite{FL83}, which characterizes the positive polynomials for ample vector bundles of rank $ r $ as those belonging to the positive convex cone spanned by the \emph{Schur polynomials} (see Section~\ref{sect: generalized schur classes/forms} for a definition).
	
	Taking a vector bundle $ (E,h) \to X $ as before, we can formally evaluate the Schur polynomials on the Chern forms of $ (E,h) $, obtaining a family of characteristic differential forms on $ X $.
	These are the so-called \emph{Schur forms} of $ (E,h) $.
	The Schur forms are special representatives of the \emph{Schur classes} of $ E $, which are the cohomology classes of $ X $ obtained by formally evaluating the Schur polynomials in $ c_1(E),\dots,c_r(E) $.
	Significant examples are the Chern classes themselves, classes of the form $ c_1(E)c_{k-1}(E) - c_k(E) $ for $ 2 \le k \le r $, and the signed Segre classes.
	Recall that the total Segre class $ 1 + s_1(E) + \dots + s_n(E) $ is the inverse of the total Chern class in $ H^{\bullet}(X) $.
	
	The pointwise Hermitianized version of the above mentioned Fulton--Lazarsfeld theorem is known in the literature as \emph{Griffiths' conjecture}, which states that the positive convex cone spanned by {Schur} forms of a Griffiths positive vector bundle consists of \emph{positive} differential forms (see Section~\ref{sect: positivity of differential forms} for definitions of the three main notions of positivity for differential forms).
	
	In this context, it is straightforward to see that (any power of) the first Chern form is strongly positive for Griffiths semipositive vector bundles.
	Except for $ c_1(E,h) $, very little is known about the positivity of the other Schur forms.
	Griffiths proved in \cite[Appendix to~\S 5.(b)]{Gri69} that the second Chern form of a rank $ 2 $ Griffiths {positive} vector bundle is positive.
	However, see \cite[p.~247]{Gri69}, Griffiths himself deemed difficult to adapt his proof to the general case.
	
	In the last few years there has been a renewed interest around Griffiths' conjecture, as important partial results have appeared in the literature.
	Guler proved in \cite[Theorem~1.1]{Gul12} the positivity of the signed Segre forms, which are the Schur forms representing the above mentioned signed Segre classes.
	More recently, it has been shown in \cite[Main Application]{DF20} the \emph{strong} positivity of an entire sub-cone (which includes, for instance, all the signed Segre forms) of the Schur cone.
	Such sub-cone is spanned by all the possible push-forwards of powers of the Chern curvatures of certain natural line bundles.
	These are constructed on the (possibly incomplete) flag bundles associated with $ (E,h) $.
	For definitions and basic notions on push-forwards of differential forms see \cite[\S I]{Dem01}.
	
	Strengthening the hypotheses with stronger notions of positivity, such as (dual) Nakano semipositivity (for instance, see \cite{Fin20} for the definitions), the conjecture is proven in \cite[Proposition~3.1]{Li20} and in \cite[Theorem~1.1]{Fin20}.
	It is also worth mentioning that Griffiths himself proved his own conjecture in \cite[Proof of Theorem~D]{Gri69} for globally generated vector bundles.
	However, the problem still remains open in the general case.
	In Section~\ref{sect: comparisons between positivity} we provide an overview on the state of the art of Griffiths' conjecture, focusing, in particular, on the recent progress.
	Primarily, the purpose of such overview is to clarify which types of positivity (among those in Definition~\ref{def: positivity notions for forms}) have been obtained for the Schur forms of positive vector bundles.
	Indeed, the terminology used in the literature is currently not always standard, and sometimes may lead to some confusion.
	
	\medskip
	In this paper, after recalling the basics on positivity of forms and some needed results on push-forwards from flag bundles, we show how to obtain the positivity of the second Chern form $ c_2(E,h) $ in any rank, if $ (E,h) $ is Griffiths positive. This can be easily obtained by adapting the above mentioned Griffiths' result on the positivity of $ c_2(E,h) $ in rank $ 2 $.
	It seems that this was first observed by Guler in \cite{Gul06}, but just as a statement without proof.
	Even if the argument is simple, we prefer to give full details for the sake of completeness.
	Furthermore, we extend such result to the semipositive context in Theorem~\ref{thmpt: c2 > 0} below.
	
	Subsequently, we focus on the Griffiths' conjecture in rank $ 3 $, for which we prove the positivity of the Schur form
	\begin{equation}\label{eq:star}
	c_1(E,h) \wedge c_2(E,h) - c_3(E,h) . \tag{$\star$}
	\end{equation}
	The starting point to obtain our result is to consider \cite[Main Theorem]{DF20} in order to get the form~\eqref{eq:star} as a push-forward from the complete flag bundle of $ E $.
	But, at the best of our knowledge, we are not able to deduce the positivity of this form instantly from that push-forward.
	
	Hence, the idea presented in this paper is to factorize such push-forward through the projective bundle $ \PP(E) $ of lines in $ E $.
	Such idea comes from an alternative proof of the positivity of $ c_2(E,h) $ in rank $ 3 $, given by the author in an earlier version of this paper.
	Coupling this factorization to an explicit push-forward formula for differential forms, which we establish in Formula~\eqref{eq: darondeau-pragacz forms schur} below, we find that~\eqref{eq:star} can be obtained also by pushing forward a positive form on $ \PP(E) $, leading us to state our main result.
	\begin{thm}\label{thm: c_2 > 0}
		Let $ (E,h) \to X $ be a Griffiths semipositive Hermitian holomorphic vector bundle of rank $ 3 $ over a complex manifold.
		Then the Schur form $ c_1(E,h) \wedge c_2(E,h) - c_3(E,h) $ is a positive $ (3,3) $-form on $ X $.
	\end{thm}
	Consequently, we find that the following chain of inequalities
	\begin{equation}\label{eq:asterisk}
		{c_1(E,h)}^{\wedge 3} \ge c_1(E,h) \wedge c_2(E,h) \ge c_3(E,h) \tag{$*$}
	\end{equation}
	holds for Griffiths semipositive vector bundles of rank $ 3 $ over complex manifolds.
	Similar inequalities have appeared, for instance, in \cite[Theorem~3.2]{Li20} if $ h $ is dual Nakano semipositive, and they give a metric counterpart of some well-known inequalities between Chern numbers valid for nef vector bundles on compact K\"ahler manifolds (see, for instance, \cite{FL83,DPS94,LZ20}).
	In the setting of Theorem~\ref{thm: c_2 > 0}, if $ X $ is compact (not necessarily K\"ahler) of dimension $ 3 $, it follows from~\eqref{eq:asterisk} that the Chern numbers of $ E $ satisfy the relation
	\begin{equation*}
		{c_1}^{3} \ge c_1  c_2 \ge c_3 .
	\end{equation*}
	
	Lastly, in Section~\ref{sect: final section} we collect some concluding remarks and open questions which may serve for further developments of this topic.
	
	\subsubsection*{Acknowledgements} The author would like to thank his advisor, Simone Diverio, for the help provided during this work, and for the valuable comments on the topics of this paper.
	
	The author is member of GNSAGA (Gruppo Nazionale per le Strutture Algebriche, Geometriche e le loro Applicazioni) of the Istituto Nazionale di Alta Matematica (INdAM) and is partially supported by \textquotedblleft{Progetti per Avvio alla Ricerca 2020}\textquotedblright, Project ID: AR120172A92A1FE7, granted by Sapienza Universit\`a di Roma.

	\section{Preliminaries and notation}
	This section is dedicated to introducing the terminology we use throughout this paper.
	Furthermore, we recall the tools (essentially taken from \cite{DP17} and \cite{DF20}) serving to state Formula~\eqref{eq: darondeau-pragacz forms schur}: a push-forward formula for differential forms on (possibly incomplete) flag bundles.
	Such formula is particularly useful in the proof of Theorem~\ref{thm: c_2 > 0}.
	
	\subsection{Positivity notions for differential forms}\label{sect: positivity of differential forms}
	We recall here some basic notions about the positivity of differential forms.
	The following exposition is taken, mainly, from \cite{HK74} and \cite[\S III]{Dem01}.
	
	Let $ V $ be a complex vector space of dimension $ n $ and let $ (e_1,\dots,e_n) $ be a basis of $ V $.
	The notation $ (e_1^{\vee},\dots,e_n^{\vee}) $ stands for the corresponding dual basis of $ V^{\vee} $.
	For $ 0 \le p,q \le n $, denote by $ \Lambda^{p,q} V^{\vee} $ the space of exterior forms of bi-degree $ (p,q) $ on $ V $.
	A form $ u \in \Lambda^{p,q} V^{\vee} $ is \emph{real} if $ \bar{u} =u $; since $ \bar{u} \in \Lambda^{q,p} V^{\vee} $, it follows that $ p=q $ in this case.
	Let $ \Lambda_{\R}^{p,p} V^{\vee} $ denote the space of real $ (p,p) $-forms.
	\begin{defn}\label{def: positivity notions for volume forms}
		A form $ \nu \in \Lambda^{n,n} V^{\vee} $ is called a \emph{positive} volume form if $ \nu = \tau \ i e_1^{\vee} \wedge \overline{e}_1^{\vee} \wedge \dots \wedge i e_n^{\vee} \wedge \overline{e}_n^{\vee} $ for some $ \tau \in \R $, $ \tau \ge 0 $.
		
		Of course, this notion is independent on the choice of the basis $ (e_1,\dots,e_n) $.
	\end{defn}
	From now on, for $ 0 \le p \le n $ set $ q = n - p  $.
	Recall that a $ (q,0) $-form which can be expressed as $ \beta_1 \wedge \dots \wedge \beta_q $ for $ \beta_1,\ldots, \beta_q  \in V^{\vee} $ is called \emph{decomposable}.
	\begin{defn}\label{def: positivity notions for forms}
		An exterior form $ u \in \Lambda_{\R}^{p,p} V^{\vee} $ is called
		\begin{itemize}
			\item \emph{positive}, if for every $ \beta \in \Lambda^{q,0} V^{\vee} $ decomposable, $ u \wedge i^{q^2} \beta \wedge \bar{\beta} $ is a positive volume form;
			\item \emph{Hermitian positive}, if for every $ \beta \in \Lambda^{q,0} V^{\vee} $, $ u \wedge i^{q^2} \beta \wedge \bar{\beta} $ is a positive volume form;
			\item \emph{strongly positive}, if there are decomposable forms $ \alpha_1,\dots,\alpha_N \in \Lambda^{p,0} V^{\vee} $ such that $ u $ can be expressed as $ \sum_{s=1}^{N} i^{p^2} \alpha_s \wedge \bar{\alpha}_s $.
		\end{itemize}
	\end{defn}
	\begin{rem}
		The terminology used for the positivity of forms is, currently, not always standard in the literature.
		We have chosen here to follow the terminology of \cite[\S III]{Dem01}, which is a standard reference for positivity of forms.
		Note however that \cite[\S III]{Dem01} does not deal with the intermediate notion of positivity (which is called simply \emph{positivity} in \cite{HK74}).
		The reason why we call it \emph{Hermitian} is due to the fact that $ u \in \Lambda_{\R}^{p,p} V^{\vee} $ is Hermitian positive if and only if for each $ \beta,\eta \in \Lambda^{q,0} V^{\vee} $, $ (\beta,\eta) \mapsto {u \wedge i^{q^2} \beta \wedge \bar{\eta}} $ gives a positive semidefinite Hermitian form on $ \Lambda^{q,0} V^{\vee} $.
	\end{rem}
	\begin{rem}\label{rem: we do not use strict}
		Of course, with an appropriate modification of Definition~\ref{def: positivity notions for forms}, one can define the \emph{strict} notions of positivity for forms.
		Mainly, we are dealing here with Griffiths semipositive vector bundles, so it is more natural to consider the three notions given in Definition~\ref{def: positivity notions for forms}, rather than the strict ones.
		Moreover, note that some authors (see, for instance, \cite{Zhe00,Gul12,Li20}) use the term \emph{non-negative}, or \emph{semi-positive}, (resp. \emph{positive}) instead of \emph{positive} (resp. \emph{strictly positive}).
	\end{rem}
	\begin{ex}\label{ex: Hermitian positive standard form}
		Let $ \xi \in \Lambda^{p,0} V^{\vee} $, then $ i^{p^2} \xi \wedge \bar{\xi} $ is a Hermitian positive form.
		Indeed, for every $ \beta \in \Lambda^{q,0} V^{\vee} $, the wedge product $ \xi \wedge \beta $ equals $ \lambda e_1^{\vee} \wedge \dots \wedge e_n^{\vee} $ for some $ \lambda \in \C $, hence $ i^{p^2} \xi \wedge \bar{\xi} \wedge i^{q^2} \beta \wedge \bar{\beta} = i^{n^2} \xi \wedge \beta \wedge \overline{ {\xi} \wedge {\beta}} $ is a positive volume form.
	\end{ex}
	
	Let $ \operatorname{WP}^p V^{\vee}$, $ \operatorname{HP}^p V^{\vee} $ and $ \operatorname{SP}^p V^{\vee} $ denote respectively the closed positive convex cones contained in $ \Lambda_{\R}^{p,p} V^{\vee} $ spanned by positive, Hermitian positive and strongly positive forms.
	The notations $ \operatorname{WP} $ and $ \operatorname{SP} $ are taken from \cite{HK74}.
	It is straightforward to see that, in general,
	\begin{equation}\label{eq: inclusions positive cones}
		\operatorname{SP}^p V^{\vee} \ \subseteq \ \operatorname{HP}^p V^{\vee} \ \subseteq \ \operatorname{WP}^p V^{\vee}.
	\end{equation}
	\begin{rem}\label{rem: positivity notions coincide when dimension}
		The two inclusions in Formula~\eqref{eq: inclusions positive cones} become equalities when $ p = 0, 1, n-1, n $.
		Indeed, if $ p = 0, n $ all the positivity notions in Definition~\ref{def: positivity notions for forms} do coincide.
		If $ p = 1 $, the one-to-one correspondence between Hermitian forms and real (1, 1)-forms shows by a diagonalization argument that every positive $ (1,1) $-form is strongly positive.
		By duality, it is also true if $ p = n-1 $.
	\end{rem}
	If $ \operatorname{K} $ is a convex cone in $ \Lambda_{\R}^{q,q} V^{\vee} $, then its \emph{dual cone} is
	\[
	\operatorname{K}^{*} := \{ u \in \Lambda_{\R}^{p,p} V^{\vee} \mid u \wedge v \text{ is a positive volume form } \forall v \in \operatorname{K} \} .
	\]
	By Definition~\ref{def: positivity notions for forms}, $ \operatorname{WP}^p V^{\vee} = \left( \operatorname{SP}^q V^{\vee} \right)^{*} $ and given that the bidual of a convex cone is equal to its closure, we have that $ \operatorname{SP}^p V^{\vee} = \left( \operatorname{WP}^q V^{\vee} \right)^{*} $.
	
	We need the following characterization of Hermitian positivity, which follows from the more general \cite[Theorem~1.2]{HK74}.
	\begin{prop}\label{prop: canonical form of forms}
		A form $ u \in \Lambda_{\R}^{p,p} V^{\vee} $ is Hermitian positive if and only if there are $ \xi_1,\dots,\xi_N \in \Lambda^{p,0} V^{\vee} $ such that $ u = \sum_{s=1}^{N} i^{p^2} \xi_s \wedge \bar{\xi}_s $.
	\end{prop}
	From Proposition~\ref{prop: canonical form of forms}, we deduce that the wedge product of two Hermitian positive forms is again Hermitian positive (of course, an analogous property holds for strongly positive forms by definition).
	Moreover, it is now clear that the dual cone $ \left( \operatorname{HP}^q V^{\vee} \right)^{*} $ equals $ \operatorname{HP}^p V^{\vee} $.
	
	We also want to recall that, sometimes, Proposition~\ref{prop: canonical form of forms} is given in the literature as a definition; see for instance \cite[p.~240]{Gri69} and \cite[\S 2]{Li20}.
	
	\begin{rem}\label{rem: positivity notions differs when dimension}
		If $ 2 \le p \le n-2 $, the two inclusions in Formula~\eqref{eq: inclusions positive cones} are strict.
		To see this, it is sufficient to observe that the Hermitian positive form $ i^{p^2} \xi \wedge \bar{\xi} $ is strongly positive if and only if $ \xi $ is decomposable (see \cite[\S III,~(1.10) Remark]{Dem01} and \cite[Proposition~1.5]{HK74}).
		Given that, for instance, the $ (p,0) $-form $ (e_1^{\vee} \wedge e_2^{\vee} + e_3^{\vee} \wedge e_4^{\vee})\wedge e_5^{\vee} \wedge \dots \wedge e_{p+2}^{\vee} $ is not decomposable, we deduce that $ \operatorname{SP}^p V^{\vee} \ \subsetneq \ \operatorname{HP}^p V^{\vee} $.
		Since the duality of cones reverses the inclusions, we get also that $ \operatorname{HP}^p V^{\vee} \subsetneq \operatorname{WP}^p V^{\vee} $.
	\end{rem}
	An explicit example of a positive form which is not Hermitian positive can be found in \cite[p.~50]{HK74}.
	As a byproduct, \cite{HK74} constructs a positive form and a Hermitian positive form, for which their wedge product is a negative volume form.
	Therefore, unlike $ \operatorname{SP}^p V^{\vee} $ and $ \operatorname{HP}^p V^{\vee} $, the cone of positive forms is not stable under wedge product (compare this with \cite{BP13}).
	
	Finally, it is useful to recall the following characterization of positivity.
	\begin{prop}\label{prop: iff weakly positive}
		A form $ u \in \Lambda_{\R}^{p,p} V^{\vee} $ is {positive} if and only if, equivalently,
		\begin{enumerate}
			\item for every vector subspace $ S \subseteq V $ with $ \dim_{\C}S = p $, the restriction $ u|_S $ is a positive volume form on $ S $;
			\item for every $ w_1,\dots,w_p \in V $, $ (-i)^{p^2} u(w_1,\dots,w_p,\overline{w}_1,\dots,\overline{w}_p) \ge 0 $.
		\end{enumerate}
	\end{prop}
	
	Of course, all of the definitions and results given in this section can be extended to a complex manifold $ X $.
	It is sufficient to take $ V = T_{X,x} $ for every $ x \in X $ and to check that all the concepts expressed here are independent by change of holomorphic coordinates.
	This follows from Definition~\ref{def: positivity notions for volume forms}.
	
	If $ E \to X $ is a complex vector bundle and $ 0 \le p,q \le n $, then $ \mathcal{A}^{p,q}(X,E) $ stands for the space $ C^{\infty}\bigl(X,\Lambda^{p,q}T_X^{\vee}\otimes E\bigr) $ of differential $ (p,q) $-forms on $ X $ with values in $ E $.
	In particular, $ \mathcal{A}^{p,q}(X) $ denotes the space of differential $ (p,q) $-forms on $ X $.
	Similarly, we use the notations $ \mathcal{A}_{\R}^{p,q}(X) $ and $ \mathcal{A}^{k}(X) $.
	
	\subsection{Generalized Schur classes/forms}\label{sect: generalized schur classes/forms}
	Fix $ k \in \N $ and denote by $ \Lambda(k,r) $ the set of partitions of $ k $ in $ r $ parts, \textsl{i.e.}, the set of all $ \sigma = (\sigma_1,\ldots,\sigma_k) \in \N^k $ such that $ r \ge \sigma_1 \ge \dots \ge \sigma_k \ge 0 $ and $ |\sigma| = \sum_{j=1}^{k} \sigma_j = k $.
	To $ \sigma \in \Lambda(k,r) $ we associate the \emph{Schur polynomial} $ S_{\sigma} \in \Z[c_1,\ldots,c_r] $ of weighted degree $ 2k $ (we regard $ c_j $ as having degree $ 2j $), which is defined as
	\[
	S_{\sigma}(c_1,\ldots,c_r) := \det \bigl( c_{\sigma_i + j - i} \bigr)_{1 \le i,j \le k} .
	\]
	By convention, $ c_0 = 1 $ and $ c_{\ell} = 0 $ if $ \ell < 0 $ or $ \ell > r $.
	
	As $\sigma\in\Lambda(k,r)$ varies, the Schur polynomials form a basis for the $ \Q $-vector space of degree $ 2k $ weighted homogeneous polynomials in $ r $ variables, and the product of two Schur polynomials is a positive linear combination of Schur polynomials.
	Hence, Schur polynomials generate a positive convex cone closed under product.
	
	From now on, let $ X $ be a complex manifold of dimension $ n $ and let $ (E,h) $ be a Hermitian holomorphic vector bundle of rank $ r $  over $ X $.
	The definitions of Schur forms (resp. classes) of $ (E,h) $ (resp. $ E $) are already given in the Introduction.
	We denote by $ S_{\sigma}(E,h) $ the {Schur form} associated to the partition $ \sigma $.
	Similarly, $ S_{\sigma}(E) $ is the {Schur class} associated to $ \sigma $.
	Clearly, $ S_{\sigma}(E,h) $ represents the cohomology class $ S_{\sigma}(E) $.
	For instance,
	\begin{equation*}
		S_{\sigma}(E,h) =
		\begin{cases}
			c_k(E,h) &\text{if} \ \sigma = (k,0,\dots,0), \\
			c_1(E,h) \wedge c_{k-1}(E,h) - c_k(E,h) &\text{if} \ \sigma = (k-1,1,0,\dots,0), \\
			(-1)^k s_k(E,h) &\text{if} \ \sigma = (1,\dots,1),
		\end{cases}
	\end{equation*}
	where, in all the three cases above, $ \sigma \in \Lambda(k,r) $, and $ s_k(E,h) $ stands for the $ k $-th Segre form of $ (E,h) $, which represents the Segre class $ s_k(E) $.
	
	Now, we want to introduce a larger family of cohomology classes (resp. differential forms) which are given by dropping the assumption that $ \sigma $ is a partition in $ \Lambda(k,r) $.
	The following notation is taken from \cite[\S 4]{DP17}, and will be essential in order to state Formula~\eqref{eq: darondeau-pragacz forms schur}.
	\begin{defn}\label{def: generalized schur classes}
		Let $ \sigma = (\sigma_1,\ldots,\sigma_k) \in \Z^k $ be a sequence of integers.
		We define a cohomology class
		\[
		s_{\sigma}(E) := \det \bigl( s_{\sigma_i + j - i}(E) \bigr)_{1 \le i,j \le k}
		\]
		in $ H^{2 |\sigma|}(X,\Z) $, where, as usual, $ s_0(E) = 1 $ and $ s_\ell(E) = 0 $ if $ \ell \notin [0,n] $.
	\end{defn}
	The relationship between Schur classes and Definition~\ref{def: generalized schur classes} is given in the following.
	\begin{ex}
			Let $ \sigma $ be a partition in $ \Lambda(k,r) $.
			It follows from the well-known Jacobi--Trudi identities that
			\begin{equation}\label{eq: relation schur vs general schur}
			s_{\sigma}(E) = (-1)^{|\sigma|} S_{\sigma^\prime}(E),
			\end{equation}
			where $ \sigma^{\prime} $ is the conjugate partition of $ \sigma $, obtained through the transposition of the Young diagram of $ \sigma $.
			For instance, by Formula~\eqref{eq: relation schur vs general schur} the partition $ (1,\dots,1) $ gives $ (-1)^k c_k(E) $, which is the $k$-th signed Chern class, while the Segre class $ s_k(E) $ is associated to the partition $ (k,0,\dots,0) $.
	\end{ex}
	Note that the sign $ (-1)^{|\sigma|} $ in Formula~\eqref{eq: relation schur vs general schur} is due to the fact that the definition of Schur classes involves Chern classes, while Definition~\ref{def: generalized schur classes} is given in terms of Segre classes.
	The reason why we choose the Segre classes' approach in Definition~\ref{def: generalized schur classes} is to make our notation uniform with that of \cite{DP17}.
	
	Similarly to Definition~\ref{def: generalized schur classes}, we can associate to any $ \sigma \in \Z^k $ the differential form $ s_{\sigma}(E,h) := \det\bigl(s_{\sigma_i+j-i}(E,h)\bigr) $.
	Inspired by Formula~\eqref{eq: relation schur vs general schur}, we call $ s_{\sigma}(E,h) $ (resp. $ s_{\sigma}(E) $) the \emph{generalized Schur form} (resp. \emph{class}) associated to $ \sigma $.
	
	\subsection{Flag bundles}
	Now we introduce some notation about flag bundles associated to the holomorphic vector bundle $ E \to X $.
	Fix a sequence of integers $ \rho=(\rho_0,\ldots,\rho_m) $ of the form $ 0 = \rho_0 < \rho_1 < \ldots < \rho_{m-1} < \rho_m = r $, the \emph{flag bundle} of $ E $ of type $ \rho $ is the holomorphic fiber bundle
	\[
	\pi_{\rho} \colon \FF_{\rho}(E) \to X
	\]
	where the fiber over $ x \in X $ is the flag manifold $ \FF_{\rho}(E_x) $, whose points are flags of the form $ \{ 0_{x} \} = V_{x,\rho_0} \subset \cdots \subset V_{x,\rho_j} \subset \cdots \subset V_{x,\rho_m} = E_{x} $ with $ \dim_{\C} V_{x,\rho_j} = \rho_j $.
	Over $ \FF_{\rho}(E) $ we have a tautological flag
	\begin{equation}\label{eq: tautological filtration vector bundles over flag bundle}
		U_{\rho,0} \subset \cdots \subset U_{\rho,j} \subset \cdots \subset U_{\rho,m}
	\end{equation}
	of vector sub-bundles of $ \pi_{\rho}^{*}E $, where, for every $ 0 \le j \le m $, the fiber of $ U_{\rho,j} $ over the point $ \{ 0_{x} \} \subset \cdots \subset V_{x,\rho_j} \subset \cdots \subset E_{x} $ is $ V_{x,\rho_j} $, hence, the rank of $ U_{\rho,j} $ is $ \rho_j $.
	We note, in particular, that $ U_{\rho,0} = (0) $ and $ U_{\rho,m} = \pi_{\rho}^{*}E $.
	\begin{ex}\label{ex: projectivized}
		If $ \rho $ is the sequence $ (0,1,r) $ then $ \FF_{\rho}(E) $ equals the projective bundle of lines in $ E $, which we denote by $ \PP(E) $.
		In this case, the tautological filtration~\eqref{eq: tautological filtration vector bundles over flag bundle} consists of one proper sub-bundle only, namely $ U_{(0,1,r),1} $, which equals, by definition, the tautological line bundle $ \OO_{\PP(E)}(-1) $.
	\end{ex}

	In the particular case of \emph{complete} flag bundles, \textsl{i.e.} $m=r$, we shall drop the subscript $\rho$ and simply write $\pi$, $ \FF(E)$ and $ U_{j} $.
	
	We denote by $ \pi_{\rho}^{r} \colon \FF(E) \to \FF_{\rho}(E) $ the obvious projection  which sends the complete flag $\{ 0_{x} \} \subset V_{x,1} \subset \cdots \subset V_{x,r-1} \subset E_{x} $ to the (partial) flag $ \{ 0_{x} \} \subset V_{x,\rho_1} \subset \cdots \subset V_{x,\rho_{m-1}} \subset E_{x} $ of $ \FF_{\rho}(E) $.
	
	Finally, note that the natural commutative diagram
	\begin{equation}\label{eq: commutative diagram of flags}
	\begin{tikzcd}
	\FF(E) \arrow[dr, "\pi"'] \arrow[rr,"\pi_{\rho}^{r}"]{}
	& & \FF_{\rho}(E) \arrow[dl,"\pi_{\rho}"] \\
	& X
	\end{tikzcd}
	\end{equation}
	of projections between flag bundles induces, for $ 0 \le k \le n $, a commutative diagram of push-forwards (intended as integration along the fibers) at the level of differential forms (see \cite[\S I,~(2.14) Theorem]{Dem01})
	\begin{equation*}\label{eq: commutative diagram of flags differential forms}
	\begin{tikzcd}
	\mathcal{A}^{2(d + k)}\bigl(\FF(E)\bigr) \arrow[dr, "\pi_{*}"'] \arrow[rr,"(\pi_{\rho}^{r})_{*}"]{}
	& & \mathcal{A}^{2(d_{\rho} + k)}\bigl(\FF_{\rho}(E)\bigr) \arrow[dl,"(\pi_{\rho})_{*}"] \\
	& \mathcal{A}^{2k}(X)
	\end{tikzcd}
	\end{equation*}
	where $ d $ (resp. $ d_{\rho} $) stands for the relative dimension of $ \pi $ (resp. $ \pi_{\rho} $).
	We also use the same notation, namely $ \pis $ and $ (\pi_{\rho})_{*} $, for the push-forwards induced in cohomology by integration along the fibers.
	
	\subsection{Push-forward formul{\ae}}\label{sect: push-forward formulas}
	Following again the notation of \cite{DP17}, denote by $ \xi_1, \ldots, \xi_r $ the (virtual) Chern roots of $ \pi_{\rho}^{*}E^{\vee} $.
	Let $ P $ be a polynomial with the appropriate symmetries which ensure that $ P(\xi_1,\ldots,\xi_r) \in {H}^{\bullet}\bigl(\FF_{\rho}(E)\bigr) $.
	In this setting, we now recall the determinantal version of Darondeau--Pragacz formula \cite[Proposition~4.2]{DP17}, given with a little modification as follows.
	
	\begin{prop}\label{prop: push-forward classes schur}
		Take the polynomial $ P $ as before, and write the class $ P(\xi_1,\ldots,\xi_r) $ as $ \sum a_{\lambda} \xi_1^{\lambda_1} \cdots \xi_r^{\lambda_r} $.
		If $ \nu $ is the increasing sequence of integers determined by $ \rho $ as:
		\[
		\nu_i = r - \rho_s \quad \text{for } r - \rho_s < i \le r - \rho_{s-1}
		\]
		(in particular, if $ m=r $ the sequence $ \nu $ is given by $ \nu_i = i-1 $), then, in terms of generalized Schur classes, we have
		\begin{equation}\label{eq: darondeau-pragacz classes schur}
		(\pi_{\rho})_{*} P(\xi_1,\ldots,\xi_r) = \sum a_{\lambda} s_{(\lambda-\nu)^{\leftarrow}}(E)
		\end{equation}
		where the notation $ (\sigma_1,\ldots,\sigma_r)^{\leftarrow} $ stands for $ (\sigma_r,\ldots,\sigma_1) $ and the difference of $ \lambda $ and $ \nu $ is defined componentwise.
	\end{prop}
	
	\begin{rem}
		Proposition~\ref{prop: push-forward classes schur} is stated a little differently from \cite[Proposition~4.2]{DP17}. 
		Indeed, the original statement of Darondeau and Pragacz is given in terms of the Chern roots of $ U_{\rho,m-1}^{\vee} $ (which is the greatest proper tautological sub-bundle of $ \pi_{\rho}^{*}E^{\vee} $), while we need a push-forward formula involving all the Chern roots of $ \pi_{\rho}^{*}E^{\vee} $.
		
		However, we have chosen to omit the proof of Proposition~\ref{prop: push-forward classes schur} since it is essentially the same as \cite[Proposition~4.2]{DP17}.
		The only relevant difference is that we use \cite[Proposition~1.2]{DP17} in the first step of the proof instead of \cite[Theorem~1.1]{DP17}.
	\end{rem}
	
	We now move to the differential forms level and suppose thereafter that $ (E,h) \to X $ is a Hermitian holomorphic vector bundle.
	
	Consider the filtration~\eqref{eq: tautological filtration vector bundles over flag bundle} of tautological vector bundles over $ \FF_{\rho}(E) $, all of them equipped with the respective restriction metrics induced by $ h $.
	One can form, for $  j=1,\dots,m $, tautological line bundles
	$$
	Q_{\rho,j}:=\det \bigl( U_{\rho,m-j+1}/U_{\rho,m-j} \bigr)
	$$
	endowed with the determinant of the quotient metrics.
	By a slight abuse of notation, denote again by $ h $ all these above mentioned metrics, and let
	$$ \Xi_{\rho,j} := \frac{i}{2\pi} \Theta(Q_{\rho,j},h) \in \mathcal{A}_{\R}^{1,1}\bigl(\FF_{\rho}(E)\bigr) $$
	be the first Chern form $ c_1(Q_{\rho,j},h) $.
	Of course, $ \Xi_{\rho,j} $ represents $ c_1(Q_{\rho,j}) $.
	
	When $m=r$, we drop the subscript $\rho$ and simply write $ Q_j $ and $ \Xi_j $.
	Note that, in this case, $ \xi_j = - c_1(Q_j) = - \left[ \Xi_j \right] $.
	
	Given a homogeneous polynomial $F$ in $m$ variables of degree $d_\rho+k$, we now want to compute the push-forward  $ (\pi_{\rho})_{*} F(\Xi_{\rho,1},\dots,\Xi_{\rho,m}) $ in terms of generalized Schur forms of $ (E,h) $.
	Clearly, at the cohomology level we have that $ F \bigl( c_1(Q_{\rho,1}), \dots, c_1(Q_{\rho,m}) \bigr) \in {H}^{2(d_{\rho} + k)}\bigl(\FF_{\rho}(E)\bigr) $.
	Therefore, there exists a polynomial in the virtual Chern roots of $ \pi_{\rho}^{*}E^{\vee} $ called
	$$
	f(\xi_1,\ldots,\xi_r) = \sum_{|\lambda| = d_{\rho}+k} b_{\lambda} \xi_{1}^{\lambda_1} \cdots \xi_{r}^{\lambda_r}
	$$
	such that, for $ 0 \le j \le m $ and $ s_j := r - \rho_{m-j} $,
	\begin{equation*}\label{rem: symmetries of F vs f} 
	\begin{aligned}
	f(\xi_1,\ldots,\xi_r) = F\left(\ldots,-\sum_{\ell = s_{j-1}+1}^{s_j} \xi_\ell,\ldots\right) = F \bigl( c_1(Q_{\rho,1}), \dots, c_1(Q_{\rho,m}) \bigr) ,
	\end{aligned}
	\end{equation*}
	\textsl{i.e.}, the polynomial $ f $ has by construction the appropriate symmetries for which $ f(\xi_1,\ldots,\xi_r) $ can be considered as a cohomology class of $ \FF_{\rho}(E) $.
	
	Note in particular that, if $ m=r $, the Chern roots of $ \pi^{*}E^{\vee} $ are not virtual.
	Hence, the only relevant symmetry in this case is $ f = (-1)^{d_{\rho} + k} F $.
	
	Thanks to \cite[Theorem~3.5]{DF20} we are able to identify the push-forward of the form $ F(\Xi_{\rho,1},\dots,\Xi_{\rho,m}) $ on the manifold $ X $.
	Accordingly to \cite[Proposition~1.2]{DP17}, the explicit expression of such push-forward is given in \cite{DF20} as the coefficient of certain monomial of a polynomial depending on $ n $, $ r $ and $ f $, formally evaluated on the Segre forms $ s_j(E,h) $'s.
	Since the explicit expression of the push-forward is obviously independent of the method we use to compute it, we prefer for our purposes to couple \cite[Theorem~3.5]{DF20} with Proposition~\ref{prop: push-forward classes schur}, obtaining the following identity
	\begin{equation}\label{eq: darondeau-pragacz forms schur}
	(\pi_{\rho})_{*} F(\Xi_{\rho,1},\dots,\Xi_{\rho,m}) = \sum_{|\lambda| = d_\rho + k} b_{\lambda} s_{(\lambda-\nu)^{\leftarrow}}(E,h) ,
	\end{equation}
	where $ \nu $ is given as in Proposition~\ref{prop: push-forward classes schur} and the $ b_{\lambda} $'s are the coefficients of $ f $.

	\section{Proof of the main results}
	In what follows, we proceed to observe that the positivity of the second Chern form is implied by the Griffiths \emph{semi}positivity of the curvature tensor.
	Such argument is central in the proof Theorem~\ref{thm: c_2 > 0}.
	
	Let $ (E,h) $ be a Griffiths positive Hermitian holomorphic vector bundle of rank $ 2 $ over a complex manifold $ X $.
	Thanks to \cite[Appendix to~\S 5.(b)]{Gri69}, we know that the second Chern form $ c_2(E,h) $ is positive.
	The argument used by Griffiths in the proof consists in the following.
	First, by Proposition~\ref{prop: iff weakly positive}, one can assume that $ X $ is a complex surface.
	Of course, this assumption is not allowed if one wants to show the Hermitian (or strong) positivity of $ c_2(E,h) $ (cf. with Remarks~\ref{rem: positivity notions coincide when dimension} and~\ref{rem: positivity notions differs when dimension}).
	However, in rank $ 2 $ the curvature is a $ 2 \times 2 $ matrix of $ (1,1) $-forms, thus
	\[
	c_2(E,h) = - \frac{1}{4 \pi^2} \det\Theta(E,h) .
	\]
	Given that the Hermitian forms associated to the diagonal entries of the matrix $ i \Theta(E,h) $ are positive definite, one can perform a simultaneous diagonalization of such Hermitian forms.
	After the diagonalization, using that $ \dim X = 2 $, to show that $ - \det\Theta(E,h) $ is positive it is sufficient to apply the Schwarz inequality coupled with the definition of Griffiths positivity.
	
	Now, suppose that $ (E,h) $ is Griffiths semipositive.
	Once fixed a strictly positive $ (1,1) $-form $ \omega $ on $ X $ and $ \varepsilon > 0 $, we can apply the above argument to the $ \operatorname{Herm}(E) $-valued $ (1,1) $-form
	\begin{equation*}
	i\Theta(E,h) + \varepsilon \omega \otimes \operatorname{Id}_{E} .
	\end{equation*}
	Passing to the limit for $ \varepsilon \to 0 $ we get the following.
	\begin{lem}\label{lem: c_2 is positive for semipositive bundles}
		Given $ (E,h) \to X $ Griffiths semipositive of rank $ 2 $, the differential form $ c_2(E,h) $ is positive.
	\end{lem}
	
	\subsection{Positivity of $ c_2 $ in any rank}\label{sect: positivity of c2 any rank}
	Let $ (E,h) \to X $ be a Hermitian holomorphic vector bundle of rank $ r $, and denote by $ \Theta_{\alpha\beta} $ the entries of the curvature matrix.
	In general, we have that
	\begin{equation}\label{eq: c2 in summands of rank 2}
		\begin{split}
			c_2(E,h) &= \operatorname{tr}_{\operatorname{End}({\Lambda}^2 E)}\left( {\Lambda}^2 \frac{i}{2\pi} \Theta(E,h) \right) \\
			&=- \frac{1}{4 \pi^2} \sum_{1 \le \alpha < \beta \le r} \bigl( \Theta_{\alpha\alpha} \wedge \Theta_{\beta\beta} - \Theta_{\alpha\beta} \wedge \Theta_{\beta\alpha} \bigr) .
		\end{split}
	\end{equation}
	If $ (E,h) $ is Griffiths positive, every diagonal entry $ \Theta_{\alpha\alpha} $ of $ \Theta(E,h) $ gives a positive definite Hermitian form.
	Thanks to a suggestion of an anonymous referee we realized that it is possible to apply \cite[Appendix to~\S 5.(b)]{Gri69} to all the summands in the last member of Formula~\eqref{eq: c2 in summands of rank 2}.
	This because each summand is by definition the determinant of a $ 2 \times 2 $ principal sub-matrix of the curvature matrix.
	From this follows the positivity of $ c_2(E,h) $ in any rank.
	
	In the semipositive case, the positivity of $ c_2(E,h) $ is obtained by applying Lemma~\ref{lem: c_2 is positive for semipositive bundles} (instead of the above argument) to all the summands in the last member of Formula~\eqref{eq: c2 in summands of rank 2}.
	This proves the following.
	\begin{thm}\label{thmpt: c2 > 0}
		Let $ (E,h) $ be a Griffiths semipositive Hermitian holomorphic vector bundle over a complex manifold.
		Then the second Chern form $ c_2(E,h) $ is a positive $ (2,2) $-form.
	\end{thm}
	
	\subsection{Positivity of $ c_1c_2 - c_3 $ in rank $3$}
	From now on, suppose that $ r = 3 $.
	Let $ p := \pi_{(0,1,3)} \colon \PP(E) \to X $ be the projective bundle of lines in $ E $ with associated tautological short exact sequence
	\begin{equation*}
		0 \to \OO_{\PP(E)}(-1) \hookrightarrow p^{*}E \twoheadrightarrow Q := p^{*}E / \OO_{\PP(E)}(-1)  \to 0
	\end{equation*}
	over $ \PP(E) $.
	Note that the quotient bundle $ Q \to \PP(E) $ is Griffiths semipositive of rank $ 2 $ with respect to the natural quotient metric induced by $ (E,h) $ (and denoted again by $ h $).
	Consequently, the $ (1,1) $-form $ c_1(Q,h) $ is strongly positive and, by Lemma~\ref{lem: c_2 is positive for semipositive bundles}, $ c_2(Q,h) $ is positive.

	We are now ready to prove the positivity of the Schur form
	\begin{equation*}
	S_{(2,1,0)}(E,h) = c_1(E,h) \wedge c_2(E,h) - c_3(E,h) .
	\end{equation*}
	
	\begin{rem}
		The main idea behind the next proof was developed for the first time by the author in an alternative proof of the positivity of $ c_2(E,h) $ in rank $ 3 $.
		Since the positivity of $ c_2(E,h) $ for Griffiths semipositive vector bundles is the content of Theorem~\ref{thmpt: c2 > 0}, we omit such an alternative proof to make te exposition more compact.
	\end{rem}
	
	\begin{proof}[Proof of Theorem~\ref{thm: c_2 > 0}]
		Following the notation of Section~\ref{sect: push-forward formulas}, let $ \xi_1,\xi_2,\xi_3 $ be the Chern roots of $ \pi^{*}E^{\vee} $.
		They are given by the tautological filtration~\eqref{eq: tautological filtration vector bundles over flag bundle} (where $ r=m=3 $), for which $ \pi^* E $ splits (non canonically) as a differentiable vector bundle over $ \FF(E) $ as $ Q_1 \oplus Q_2 \oplus Q_3 $.
		
		By the commutativity of the following diagram of projections
		\begin{equation}\label{eq: commutative diagram of flags in rank 3}
			\begin{tikzcd}
				\FF(E) \arrow[dr, "\pi"'] \arrow[rr,"q := \pi_{(0,1,3)}^{3}"]{}
				& & \PP(E) \arrow[dl,"p"] \\
				& X
			\end{tikzcd}
		\end{equation}
		(cf. with the commutative diagram~\eqref{eq: commutative diagram of flags} above) we see that $ \xi_1 $ and $ \xi_2 $ are also the Chern roots of $ q^{*}Q^{\vee} $, where $ Q \to \PP(E) $ is the rank $ 2 $ tautological quotient bundle introduced before.
		Indeed,
		\begin{equation*}
			q^* Q = q^* \left( p^*E / \OO_{\PP(E)}(-1) \right) = \pi^*E/U_1
		\end{equation*}
		and, since $ \FF(E) $ coincides with $ \FF(Q) = \PP(Q) $, we deduce that there is a non canonical isomorphism
		\[
		q^* Q \cong_{C^{\infty}} \pi^*E/U_2 \oplus U_2/U_1 = Q_1 \oplus Q_2.
		\]
		
		By Formula~\eqref{eq: darondeau-pragacz forms schur} applied to $ \pi \colon \FF(E) \to X $ and where $ F(\Xi_1,\Xi_2,\Xi_3) $ is taken to be the monomial $ \Xi_1^{\wedge 4} \wedge \Xi_2^{\wedge 2} \wedge \Xi_3^{\wedge 0} $, we have the equality
		\begin{equation}\label{eq: c2 at form level, flag}
			\pis[\Xi_1^{\wedge 4} \wedge \Xi_2^{\wedge 2} \wedge \Xi_3^{\wedge 0}] = s_{(-2,1,4)}(E,h) = S_{(2,1,0)}(E,h) .
		\end{equation}
		Observe however that, \textsl{a priori}, we do not conclude anything about the positivity of $ S_{(2,1,0)}(E,h) $ using only Equation~\eqref{eq: c2 at form level, flag}, since, for instance, the line bundle $ Q_2 $ is not positive in general.
		Therefore, we apply again Formula~\eqref{eq: darondeau-pragacz forms schur}, but to the flag bundle $ q \colon \PP(Q) \to \PP(E) $, getting
		\begin{equation}\label{eq: square of the second chern form of quotient}
			q_{*}[\Xi_1^{\wedge 4} \wedge \Xi_2^{\wedge 2}] = s_{(1,4)}(Q,h) = c_1(Q,h) \wedge c_2(Q,h)^{\wedge 2} ,
		\end{equation}
		where the last equality holds given that $ c_3(Q,h) \equiv 0 $ and $ c_4(Q,h) \equiv 0 $.
		
		Hence, the chain of equalities
		\begin{align*}
			p_{*} \left[ c_1(Q,h) \wedge c_2(Q,h)^{\wedge 2} \right] &= p_{*} q_{*}[\Xi_1^{\wedge 4} \wedge \Xi_2^{\wedge 2}] &\\
			&= \pis [\Xi_1^{\wedge 4} \wedge \Xi_2^{\wedge 2} \wedge \Xi_3^{\wedge 0}] &\text{by commutativity of~\eqref{eq: commutative diagram of flags in rank 3}}\\
			&= S_{(2,1,0)}(E,h) &\text{by Formula~\eqref{eq: c2 at form level, flag}}
		\end{align*}
		follows, by applying $ p_{*} $ to both members of Formula~\eqref{eq: square of the second chern form of quotient}.
		
		Since $ Q \to \PP(E) $ is Griffiths semipositive, by Lemma~\ref{lem: c_2 is positive for semipositive bundles} we know that $ c_2(Q,h) $ is a positive form on $ \PP(E) $.
		Thanks to \cite[Theorem~1]{BP13}, the square of a positive $ (2,2) $-form is also positive.
		Hence $ c_1(Q,h) \wedge c_2(Q,h)^{\wedge 2} $ is positive, given that it is the wedge product of a strongly positive and of a positive form.
		Consequently, $ S_{(2,1,0)}(E,h) $ must be a positive form on $ X $, since it is the push-forward of a positive differential form.
	\end{proof}
	
	\begin{rem}\label{rem: when c2 is strongly positive}
		For dimension and bi-degree reasons (see Remark~\ref{rem: positivity notions coincide when dimension}) we underline that if $ \dim X \le 3 $, then $ c_2(E,h) $ is strongly positive, and  if $ \dim X \le 4 $, then $ S_{(2,1,0)}(E,h) $ is strongly positive, too.
	\end{rem}
	
	From the (strong) positivity of the second Segre form
	\[
	s_2(E,h) = c_1(E,h)^{\wedge 2} - c_2(E,h)
	\]
	(see \cite{Gul12,DF20}) and from Theorem~\ref{thm: c_2 > 0} we deduce the following.
	\begin{cor}\label{cor: chain of equalities}
		If $ (E,h) \to X $ is Griffiths semipositive of rank $ 3 $ over a complex manifold, then the following pointwise inequalities hold
		\begin{equation}\label{eq: inequalities chern forms as corollary}
			{c_1(E,h)}^{\wedge 3} \ge c_1(E,h) \wedge c_2(E,h) \ge c_3(E,h) .
		\end{equation}
		In addition, if $ X $ is compact and $ 3 $-dimensional, the following chain of inequalities between Chern numbers also holds
		\begin{equation}\label{eq: inequalities chern numbers as corollary}
			\int_{X} {c_1(E)}^{3} \ge \int_{X} c_1(E) c_2(E) \ge \int_{X} c_3(E) .
		\end{equation}
	\end{cor}
	As already mentioned in the Introduction, the inequalities in~\eqref{eq: inequalities chern forms as corollary} can be deduced by \cite[Theorem~3.2]{Li20} only if $ h $ is dual Nakano semipositive.
	Thus, Corollary~\ref{cor: chain of equalities} is in some sense a generalization of \cite[Theorem~3.2]{Li20} in rank~$ 3 $.
	\begin{rem}
		It is well-known that a vector bundle admitting a Griffiths semipositive metric is nef (see \cite{DPS94} for a definition).
		Recall that the Chern numbers of a nef vector bundle on a compact $ n $-dimensional K\"ahler manifold are bounded above by the Chern number $ c_1^n $ (\cite[Corollary~2.6]{DPS94}) and below by the Euler number $ c_n $ (\cite[Theorem~2.9]{LZ20}). See also \cite[Remark~3.3]{Li20}.
		Therefore, assuming that $ X $ is also K\"ahler the first and the second inequalities in~\eqref{eq: inequalities chern numbers as corollary} are a particular case of \cite[Corollary~2.6]{DPS94} and \cite[Theorem~2.9]{LZ20} respectively.
		Moreover, observe that the inequalities in~\eqref{eq: inequalities chern numbers as corollary} follow from \cite[Theorem~2.5]{DPS94} if $ X $ is K\"ahler.
	\end{rem}
	
	\begin{rem}\label{rem: difficulty of c3}
		It seems difficult to apply the strategy in the proof of Theorem~\ref{thm: c_2 > 0} to prove the positivity of $ c_3(E,h) $ in rank $ 3 $.
		Indeed, in order to factor through the quotient bundle, we need to push-forward only monomials of the form $ \Xi_1^{\wedge \lambda_1} \wedge \Xi_2^{\wedge \lambda_2} $.
		But, by Formula~\eqref{eq: darondeau-pragacz forms schur} we have the equality $ \pis[\Xi_1^{\wedge 3} \wedge \Xi_2^{\wedge 2} \wedge \Xi_3] = c_3(E,h) $ that involves $ \Xi_3 $, which is the pull-back through $ q $ of the Chern curvature of $ \OO_{\PP(E)}(-1) $.
		
		However, one can try to adapt the same ideas presented in this work to vector bundles of rank higher than $ 3 $.
		For example, in rank $ 4 $ we could use the tautological quotient of the Grassmann bundle $ \GG_2(E) $ in order to get push-forward formul{\ae} for some Schur forms.
		The problem here is that the monomials we push-forward contain $ c_2^{\wedge \ell} $, $ \ell > 2 $.
		Therefore, as pointed out for instance in \cite{BP13}, we do not know \textsl{a priori} if these powers are positive.
	\end{rem}

	\section{State of the art, concluding remarks and open questions}\label{sect: final section}
	As always, let $ (E,h) $ be a rank $ r $ Hermitian holomorphic vector bundle over a complex manifold $ X $. 
	
	\subsection{A comparison on the positivity of Schur forms}\label{sect: comparisons between positivity}
	We take the opportunity here to compare the different notions of positivity (see Definition~\ref{def: positivity notions for forms}) of Schur forms that appear in \cite{BC65,Gri69,Gul12,Li20,Fin20,DF20}, which are the works mainly related to Theorem~\ref{thm: c_2 > 0}.
	
	Recall that we do not deal with the strict notions of positivity (see Remark~\ref{rem: we do not use strict}).
	For this, the following exposition can be made slightly more precise, although this does not affect the purpose of the comparison.
	
	\subsubsection*{Positivity in \cite{BC65}}
	In \cite[Definition~5.1]{BC65} a notion of positivity for elements in $ \mathcal{A}^{p,p}\bigl(X,\operatorname{End}(E)\bigr) $ is given.
	In particular, for $ p=1 $, such notion includes the so-called Bott--Chern nonnegativity (this terminology is due to \cite{Li20}, see below), which requires that the curvature of $ (E,h) $ can be expressed, locally, as $ A \wedge \bar{A}^t $, where $ A $ is a matrix of $ (1,0) $-forms of appropriate size.
	Although \cite{BC65} predates Griffiths' conjecture, from \cite[Lemma~5.3]{BC65} one can deduce the Hermitian positivity of the top Chern form of a Bott--Chern nonnegative vector bundle.
	
	\subsubsection*{Positivity in \cite{Gri69}}
	The positivity notion considered in \cite[p.~240]{Gri69} is, by Proposition~\ref{prop: canonical form of forms}, Hermitian positivity.
	As already mentioned, Griffiths' conjecture first appears in \cite[p.~247]{Gri69}, where it is conjectured that for a Griffiths positive vector bundle the cone of positive polynomials in the Chern forms (\textsl{i.e.}, the Schur cone) consists of Hermitian positive differential forms.
	The full conjecture is verified in \cite[p.~246]{Gri69} for globally generated vector bundles; for more details we refer to \cite[Proof of Theorem~D]{Gri69}.
	
	Moreover, by the characterization of Proposition~\ref{prop: iff weakly positive}, Griffiths proves in \cite[Appendix to~\S 5.(b)]{Gri69} the positivity, but not the Hermitian positivity, of $ c_2(E,h) $, for $ (E,h) $ Griffiths positive of rank $ 2 $.
	
	\subsubsection*{Positivity in \cite{Gul12}}
	The main result of \cite{Gul12} states that the signed Segre forms of a Griffiths positive vector bundle are positive.
	Actually, \cite[Theorem~1.1]{Gul12} implicitly proves that such forms are strongly positive, even though the strong positivity is not explicitly observed therein.
	
	\subsubsection*{Positivity in \cite{Li20}}
	The two positivity notions considered by \cite{Li20} are, by Propositions~\ref{prop: iff weakly positive} and~\ref{prop: canonical form of forms},  positivity and Hermitian positivity, although in \cite[\S 2]{Li20} are called \emph{nonnegativity} and \emph{strong nonnegativity} respectively.
	In the same spirit of \cite[Proof of Theorem~D]{Gri69}, \cite[Proposition~3.1]{Li20} extends the above mentioned result on globally generated vector bundles to the larger family of Bott--Chern nonnegative (see \cite[Definition~2.1]{Li20}) vector bundles, showing that the Schur forms of these bundles are Hermitian positive.
	
	\subsubsection*{Positivity in \cite{Fin20}}
	All the three notions of Definition~\ref{def: positivity notions for forms} are addressed in \cite{Fin20}, although with the terminology of \cite{HK74}.
	The full Griffiths' conjecture is verified in \cite[Theorem~1.1]{Fin20}, where it is shown the Hermitian positivity of the Schur forms of a Nakano, or dual Nakano, positive vector bundle.
	
	It is also worth to recall that \cite{Fin20} relates, in a very interesting way, Bott--Chern's notion of positivity with the dual Nakano one.
	More precisely, \cite[Theorem~2.14]{Fin20} states that a Hermitian vector bundle $ (E,h) \to X $ is dual Nakano semipositive if and only if for every $ x \in X $ there is a vector space $ V $ and a matrix of $ (1,0) $-forms $ A \in T_{X,x}^\vee \otimes \operatorname{Hom}(V,E_x) $ such that $ \Theta(E,h)_x = A \wedge \bar{A}^t $.
	
	In order to have a similar characterization for dual Nakano positivity, we observe that it is sufficient to require, in addition, the invertibility of the operator $ \tilde{A} \in \operatorname{Hom}(T_{X,x} \otimes E_x^{\vee}, V^{\vee}) $ naturally associated to $ A $.
	Therefore, both \cite{Li20} and \cite{Fin20} obtain the full Griffiths' conjecture for dual Nakano (semi)positive bundles, but with different methods.
	
	Finally, we recall that \cite[Theorem~1.3]{Fin20} establishes the equivalence between Griffiths' conjecture and an open question concerning the so-called \emph{positive semidefinite linear preservers}, see \cite[Open problem]{Fin20}.
	
	\subsubsection*{Positivity in \cite{DF20}}
	As already mentioned in the Introduction, given a Griffiths semipositive vector bundle $ (E,h) \to X $, \cite[Main Application]{DF20} shows the strong positivity of a family of differential forms in the Schur cone of $ (E,h) $.
	Such family form a sub-cone which, in the notation of Section~\ref{sect: push-forward formulas}, is spanned by all possible wedge products of all possible push-forwards
	\begin{equation}\label{eq: diverio-fagioli push forwards}
	(\pi_{\rho})_* \bigr( a_1 \Xi_{\rho,1} + \dots + a_m \Xi_{\rho,m} \bigl)^{d_{\rho} + k}
	\end{equation}
	as $ k \ge 0 $, $ {\rho} = (\rho_0,\ldots,\rho_m) $ and the sequence of integers $ a_1 \ge \dots \ge a_m \ge 0 $ vary.

	Observe that to $ \rho = (0,r-1,r) $ corresponds the bundle $ \pi_{\rho} \colon \PP(E^{\vee}) \to X $ of hyperplanes in $ E $.
	Hence $ \Xi_{\rho,1} = c_1 \bigr( \OO_{\PP(E^{\vee})}(1),h \bigl) $ and for $ (a_1,a_2)=(1,0) $ Guler's result \cite[Theorem~1.1]{Gul12} on the signed Segre forms is recovered.
	
	\subsection{Concluding remarks}
	If we assume (dual) Nakano semipositivity, then the Schur cone consists of Hermitian positive differential forms.
	Hence, it is natural to ask the following question.
	\begin{quest}\label{quest: forms outside diverio-fagioli cone}
		Beside those as in Expression~\eqref{eq: diverio-fagioli push forwards} (and wedge products of them), are there in the Schur cone other strongly positive differential forms for (dual) Nakano semipositive vector bundles?
	\end{quest}
	Of course, the best possible result would be that all the Schur cone consists of strongly positive differential forms.
	
	\medskip
	Assume now that the vector bundle is Griffiths semipositive.
	Although it is not explicitly mentioned in Griffiths' conjecture (as stated in the literature), at this point it seems interesting to ask what is the most natural notion of positivity, among those in Definition~\ref{def: positivity notions for forms}, that we can expect to hold for the Schur forms of the vector bundle.
	Certainly, the conjecture is still open even requiring the weakest notion of positivity, although, as already said, \cite[Main Application]{DF20} shows us that the strong positivity naturally appears in this context.
	In addition, it is in a sense more natural to require the strong, or Hermitian, positivity of the Schur forms.
	Indeed, as mentioned in Section~\ref{sect: generalized schur classes/forms}, positive polynomials (\textsl{i.e.}, those belonging to the Schur cone) are stable under product and so do Hermitian and strongly positive forms; while the wedge product of two positive forms is not necessarily positive (see \cite[\S III,~(1.11) Proposition]{Dem01} and \cite{BP13}).
	
	Summing up, one may ask an analogue of Question~\ref{quest: forms outside diverio-fagioli cone} for Griffiths semipositive vector bundles, wondering if, outside the cone spanned by the push-forwards~\eqref{eq: diverio-fagioli push forwards}, there are other Hermitian or strongly positive Schur forms.
	If they were all (at least) Hermitian positive one would have an affirmative answer to the original Griffiths' conjecture as stated in \cite[p.~247]{Gri69}.
	
	However, we point out that these considerations are quite optimistic.
	For instance, beside the cases listed in Remark~\ref{rem: when c2 is strongly positive}, we do not deduce the Hermitian positivity of $ c_2(E,h) $ (resp. of $ S_{(2,1,0)}(E,h) $) by the proof of Theorem~\ref{thmpt: c2 > 0} (resp. of Theorem~\ref{thm: c_2 > 0}).
	
	\subsection{Other related questions}\label{sect: open questions}
	Let $ \mathcal{E} $ be an ample vector bundle over a projective manifold $ X $.
	
	As pointed out in \cite{Pin18,Xia20}, one can study some variants of Griffiths' conjecture.
	For instance, \cite[Conjecture~1.4]{Xia20} asks whether every Schur class $ S_{\sigma}(\mathcal{E}) $ does admit a positive representative.
	Such question fits between Fulton--Lazarsfeld theorem and Griffiths' conjecture.
	Thanks to  \cite[Theorem~A]{Xia20} we know that the answer to this question is affirmative if $ |\sigma| = \dim X -1 $.
	Another partial affirmative answer to \cite[Conjecture~1.4]{Xia20} is given in \cite[Theorem~4.8]{DF20}, where it is shown that the cohomology classes in the Schur cone coming from push-forwards as in Expression~\eqref{eq: diverio-fagioli push forwards} contain a strongly positive form.
	However, the general problem is open.
	
	Moreover, as observed in \cite[Remark~1.5]{Xia20}, given a smooth representative $ \eta \in S_{\sigma}(\mathcal{E}) $ it is not clear how to find a Hermitian metric $ h $ on $ \mathcal{E} $ such that $ \eta = S_{\sigma}(\mathcal{E},h) $.
	By virtue of these remarks, one may wonder the following.
	\begin{quest}\label{quest2}
		Given an ample vector bundle $ \mathcal{E} \to X $, is it true that for any partition $ \sigma \in \Lambda(k,r) $, there exists a Hermitian metric $ {h}_{\sigma} $ on $ \mathcal{E} $ such that $ S_{\sigma}(\mathcal{E},{h}_{\sigma}) $ is (Hermitian/strongly) positive?
	\end{quest}
	If $ \dim X = 2 $ and $ \mathcal{E} $ is semistable with respect to some polarization, \cite[Theorem~1.1]{Pin18} provides an affirmative answer to Question~\ref{quest2}, finding a Hermitian metric on $ \mathcal{E} $ whose Schur forms are positive.
	However, it is not clear if such metric is Griffiths positive: see \cite[p.~633]{Pin18}.
	Clearly, an affirmative answer to Question~\ref{quest2} would imply \cite[Conjecture~1.4]{Xia20}.
	
	We would finally remark that affirmative answers to \cite[Conjecture~1.4]{Xia20} or Question~\ref{quest2} do not imply Griffiths' conjecture, given that the latter requires the positivity of the Schur forms with respect to one \emph{given} Griffiths (semi)positive metric on the vector bundle.

\bibliography{bibliography}{}

\end{document}